# CONTROLLABILITY OF REDUCED SYSTEMS

PETRE BIRTEA[1], MIRCEA PUTA[2], AND TUDOR S. RATIU[1]

**Abstract.** Sufficient conditions for the controllability of a conservative reduced system are given. Several examples illustrating the theory are also presented.

## 1. Introduction

The phase space of classical conservative mechanical system is described usually by a symplectic manifold $(M, \omega)$. The dynamics on $M$, subject to external forces, can often be written in the form of a control system as

$$(1.1) \qquad \dot{x} = f(x) + \sum_{i=1}^{m} g_i(x) u_i, \ x \in M, \ u = (u_1, ... u_m) : \mathbb{R} \to B \subset \mathbb{R}^m,$$

where

$$f(x) = \sum_{i=1}^{n} \{x_i, H\}_\omega \frac{\partial}{\partial x_i}$$

is a Hamiltonian vector field on $M$, $\{\cdot, \cdot\}_\omega$ is the Poisson bracket given by the symplectic form $\omega$, $n = \dim M$, $B$ is a bounded set, and $u$ is measurable.

Assume that on $M$ we have a free proper action of a Lie group $G$ which preserves the symplectic form $\omega$. If all the vector fields $f$ and $g_i$ are $G$-invariant then they induce reduced dynamics on the quotient manifold $M/G$

$$(1.2) \qquad \dot{\widetilde{x}} = \widetilde{f}(\widetilde{x}) + \sum_{i=1}^{m} \widetilde{g}_i(\widetilde{x}) u_i$$

and the induced vector field $\widetilde{f}$ is also Hamiltonian on the reduced Poisson manifold $\left(M/G, \{\cdot, \cdot\}_{M/G}\right)$.

The aim of this paper is to give sufficient conditions for the controllability of the system (1.2). We will give topological conditions under which the well known sufficient conditions

(i)   $\widetilde{f}$ is weakly positively Poisson stable (WPPS)
(ii)  the Lie algebra rank condition (LARC) is satisfied

---

[1]Département de Mathématiques, École Polytechnique Fédérale de Lausanne. CH–1015 Lausanne. Switzerland.
[2]Departamentul de Matematică, Universitatea de Vest din Timişoara, Romania.
*Date*: March 8, 2002.

for the controllability of the system (1.2) are satisfied. The result will be applied to the motion of three point vortices in the plane, the three-wave interaction, and two coupled planar rigid bodies.

Poisson control systems with symmetry were studied before by Sànchez de Alvarez [30] who started a systematic investigation of the relationship of the controllability between the unreduced and reduced system. In [31] she has shown that the unreduced control system is locally weakly controllable relative to symmetries if and only if the reduced control system is locally weakly controllable. Other results concerning different aspects of the relation between the unreduced and the reduced control system can be found in Jalnapurkar and Marsden [12], [13], and Bloch, Leonard, and Marsden [6].

The present work was inspired by a very nice paper of Manikonda and Krishnaprasad [14] that studies the case $M = T^*G$, with $G$ possibly non compact, but $\mathfrak{g}^*$ with closed coadjoint orbits.

## 2. Symplectic reduction

This section quickly reviews some standard results on symplectic reduction necessary in the subsequent proofs.

Consider a $2n$-dimensional symplectic manifold $(M, \omega)$ on which there is a free proper symplectic action of a Lie group $G$. Then the orbit space $M/G$ is a smooth Poisson manifold and the projection

$$\pi : (M, \{\cdot, \cdot\}_\omega) \longrightarrow \left(M/G, \{\cdot, \cdot\}_{M/G}\right)$$

is a Poisson surjective submersion. If, in addition, the Lie group $G$ is compact then $\pi$ is a closed proper map (proofs of these statements can be found in e.g. [1], [2], [9], [14], [19]).

Suppose that the $G$-action on $M$ admits an associated $Ad^*$-equivariant momentum map $J : M \longrightarrow \mathfrak{g}^*$. The *Marsden-Weinstein reduction theorem* states that if $\mu \in \mathfrak{g}^*$ is a regular value of $J$ then the smooth quotient manifold $M_\mu := J^{-1}(\mu)/G_\mu$ is symplectic with symplectic form $\omega_\mu$ characterized by

$$\pi_\mu^* \omega_\mu = i_\mu^* \omega,$$

where $G_\mu$ denotes the isotropy subgroup of $\mu$ under the coadjoint action, $i_\mu : J^{-1}(\mu) \longrightarrow M$ is the inclusion, and $\pi_\mu : J^{-1}(\mu) \longrightarrow M_\mu$ is the projection (for a proof, see the original paper [26], or [1], [19], [28]). The symplectic manifolds $(M_\mu, \omega_\mu)$ will be called *point reduced spaces*.

These point reduced spaces $M_\mu$ can be understood in a natural way as symplectic leaves of the Poisson manifold $\left(M/G, \{\cdot, \cdot\}_{M/G}\right)$. Indeed, the smooth map $j_\mu : M_\mu \longrightarrow M/G$ naturally defined by the commutative diagram

$$\begin{array}{ccc} J^{-1}(\mu) & \xrightarrow{i_\mu} & M \\ \pi_\mu \downarrow & & \downarrow \pi \\ M_\mu & \xrightarrow{j_\mu} & M/G \end{array}$$

is a Poisson injective immersion. Moreover, the $j_\mu$-images in $M/G$ of the connected components of the symplectic manifolds $(M_\mu, \omega_\mu)$ are its symplectic leaves (see [21] or [27]).

Observe that, in general, $j_\mu$ is only an injective immersion. So the topology of the image of $j_\mu$ homeomorphic to the topology of $M_\mu$ is stronger than the subspace topology induced by the ambient space $M/G$. This image topology on $j_\mu(M_\mu)$ is called the *immersed topology*. Thus, there may be sets that are compact in the immersed topology of $j_\mu(M_\mu)$ but are not compact in the subspace topology of the same set.

The proof of the next proposition requires compactness of $G$.

**Proposition 2.1.** *Suppose that the free symplectic compact $G$-action on $(M, \omega)$ admits an $Ad^*$-equivariant momentum map $J : M \longrightarrow \mathfrak{g}^*$. Then the symplectic leaves of $\left(M/G, \{\cdot, \cdot\}_{M/G}\right)$ are closed sets.*

*Proof.* Since $J^{-1}(\mu)$ is closed in $M$ and $\pi : M \longrightarrow M/G$ is a closed map (because $G$ is compact), the set $j_\mu(M_\mu) = \pi(J^{-1}(\mu))$ is closed in the topology of $M/G$. Therefore, the connected components of $j_\mu(M_\mu)$, which are the symplectic leaves of $M/G$, are also closed in the topology of $M/G$. □

We return now to the general case with $G$ non compact. Up to now we have regarded the symplectic leaves of $\left(M/G, \{\cdot, \cdot\}_{M/G}\right)$ as the $j_\mu$-images of the connected components of $M_\mu$. However, as sets,

$$j_\mu(M_\mu) = J^{-1}(\mathcal{O}_\mu)/G,$$

where $\mathcal{O}_\mu \subset \mathfrak{g}^*$ is the coadjoint orbit through $\mu$. The set $M_{\mathcal{O}_\mu} := J^{-1}(\mathcal{O}_\mu)/G$ is called *the orbit reduced space* associated to the coadjoint orbit $\mathcal{O}_\mu$. The smooth manifold structure (and hence the topology) on $M_{\mathcal{O}_\mu}$ is the one that makes $j_\mu : M_\mu \longrightarrow M_{\mathcal{O}_\mu}$ into a diffeomorphism. The next theorem characterizes the symplectic form and the Hamiltonian dynamics on $M_{\mathcal{O}_\mu}$.

**Theorem 2.2.** *(Symplectic orbit reduction) Assume that the free proper symplectic action of the Lie group $G$ on the symplectic manifold $(M, \omega)$ admits an associated $Ad^*$-equivariant momentum map $J : M \longrightarrow \mathfrak{g}^*$.*

*(i) On $J^{-1}(\mathcal{O}_\mu)$ there is a unique immersed smooth manifold structure such that $\pi_{\mathcal{O}_\mu} : J^{-1}(\mathcal{O}_\mu) \longrightarrow M_{\mathcal{O}_\mu}$ is a surjective submersion, where $M_{\mathcal{O}_\mu}$ is endowed with the manifold structure making $j_\mu$ into a diffeomorphism. This smooth manifold structure does not depend on the choice of $\mu$ in the coadjoint orbit. If $J^{-1}(\mathcal{O}_\mu)$ is* a submanifold of $M$ in its own right then the immersed topology and the induced topology on $M_{\mathcal{O}_\mu}$ coincide.

(ii) $M_{\mathcal{O}_\mu}$ is a symplectic manifold with the symplectic form $\omega_\mu^{\mathbf{\Phi}}$ uniquely characterized by the relation

$$i_{\mathcal{O}_\mu}^* \omega = \pi_{\mathcal{O}_\mu}^* \omega_\mu^{\mathbf{\Phi}} + J_{\mathcal{O}_\mu}^* \omega_{\mathcal{O}_\mu}^+,$$

where $J_{\mathcal{O}_\mu}$ is the restriction of $J$ to $J^{-1}(\mathcal{O}_\mu)$, $i_{\mathcal{O}_\mu} : J^{-1}(\mathcal{O}_\mu) \longrightarrow M$ is the inclusion, and $\omega_{\mathcal{O}_\mu}^+$ is the $+$orbit symplectic form on $\mathcal{O}_\mu$.

(iii) Let $H$ be a $G$-invariant function on $M$ and define $\widetilde{H} : M/G \longrightarrow \mathbb{R}$ by $H = \widetilde{H} \circ \pi$. Then the Hamiltonian vector field $X_H$ is also $G$-invariant and hence induces a vector field on $M/G$ which coincides with the Hamiltonian vector field $X_{\widetilde{H}}$. Moreover, the flow of $X_{\widetilde{H}}$ leaves the symplectic leaves $M_{\mathcal{O}_\mu}$ of $M/G$ invariant. This flow restricted to the symplectic leaves is again Hamiltonian relative to the symplectic form $\omega_\mu^{\mathbf{\Phi}}$ and the Hamiltonian function $\widetilde{H}_{\mathcal{O}_\mu}$ given by

$$\widetilde{H}_{\mathcal{O}_\mu} \circ \pi_{\mathcal{O}_\mu} = H \circ i_{\mathcal{O}_\mu}.$$

The proof of this theorem in the regular case and when $\mathcal{O}_\mu$ is a submanifold of $\mathfrak{g}^*$ can be found in Marle [22], Kazhdan, Kostant, and Sternberg [15], and Marsden [23]. For the general case, when $\mathcal{O}_\mu$ is not a submanifold of $\mathfrak{g}^*$ see the book of Ortega and Ratiu [27]. Here is the main idea of the proof. Consider for each value $\mu \in \mathfrak{g}^*$ of $J$ the $G$-equivariant bijection

$$s : G \times_{G_\mu} J^{-1}(\mu) \to J^{-1}(\mathcal{O}_\mu)$$

$$[g, m] \mapsto g \cdot m$$

On $J^{-1}(\mathcal{O}_\mu)$ consider the smooth manifold structure that makes the bijection $s$ into a diffeomorphism. Then $J^{-1}(\mathcal{O}_\mu)$ with this smooth structure is an immersed submanifold of $M$.

In the particular case in which $J^{-1}(\mathcal{O}_\mu)$ is a smooth submanifold of $M$ in its own right, its smooth structure coincides with the one induced by the mapping $s$ since in this situation this bijection becomes a diffeomorphism.

If $\mu$ is a regular value of $J$ and $\mathcal{O}_\mu$ is an embedded submanifold of $\mathfrak{g}^*$ then $J$ is transversal to $\mathcal{O}_\mu$ and hence $J^{-1}(\mathcal{O}_\mu)$ is automatically an embedded submanifold of $M$.

The following statement will be used in the sequel.

**Proposition 2.3.** *(Bifurcation Lemma) Let $(M, \omega)$ be a symplectic manifold and let $G$ be a Lie group acting symplectically on $M$ (not necessarily freely). Suppose also that the action has an associated $\mathrm{Ad}^*$-equivariant momentum map $J : M \longrightarrow \mathfrak{g}^*$. For any $m \in M$,*

$$(\mathfrak{g}_m)^\circ = range(T_m J)$$

*where $\mathfrak{g}_m = \{\xi \in \mathfrak{g} \mid \xi_M(m) = 0\}$ is the Lie algebra of the isotropy subgroup $G_m = \{g \in G \mid g \cdot m = m\}$ and $(\mathfrak{g}_m)^\circ = \{\mu \in \mathfrak{g}^* \mid \mu_{|\mathfrak{g}_m} = 0\}$ denotes the annihilator of $\mathfrak{g}_m$ in $\mathfrak{g}^*$.*

An immediate consequence of this is the fact that when the action of $G$ is free then every value $\mu \in \mathfrak{g}^*$ of the momentum map $J$ is a regular value.

## 3. Controllability and Poisson stability

Let $M$ be a smooth $n$-dimensional connected manifold and

$$(3.1) \qquad \dot{x} = f(x) + \sum_{i=1}^{m} g_i(x) u_i$$

be a nonlinear control system on $M$, where $f, g_1, ..., g_m \in \mathcal{X}(M)$ are smooth vector fields on $M$, the control $u := (u_1, ..., u_m) : (0, \infty) \longrightarrow B \subset \mathbb{R}^m$ is a measurable function, and $B$ is a bounded subset of $\mathbb{R}^m$. We will denote by $\mathcal{L}$ the Lie subalgebra of $\mathcal{X}(M)$ generated by the vector fields $f, g_1, ..., g_m$, i.e.,

$$\mathcal{L} := span\,(f, g_1, ..., g_m)$$

**Definition 3.1.** *The system (3.1) satisfies the* Lie algebra rank condition (LARC) *if $\mathcal{L}(x) = T_x M$ for every $x \in M$, where $\mathcal{L}(x) := \{X_x \mid X \in \mathcal{L}\}$.*

**Definition 3.2.** *The system (3.1) is* controllable *if for any two points $x_I, x_F \in M$, there is a control $u$ which takes the system from point $x = x_I$ at time $t = t_I$ to the point $x = x_F$ at time $t = t_F$, that is, if for a certain choice of the function $u$ there is an integral curve $x(t)$ of (3.1) that begins at $x_I$ and ends at $x_F$.*

It is well known that for a nonlinear control system without drift, i.e., $f = 0$, the (LARC) condition implies controllability. This is Chow's theorem [10]. For the general case $f \neq 0$, the situation is more complicated and, in general, (LARC) is not sufficient to guarantee controllability. We will review below what is known about this case.

Let $X \in \mathcal{X}(M)$ be a smooth complete vector field on $M$ and let $\{\Phi_t\}_{t \in \mathbb{R}}$ be its flow.

**Definition 3.3.** *A point $x \in M$ is called* positively Poisson stable *for $X \in \mathcal{X}(M)$ if for any $T > 0$ and any neighborhood $V_x$ of $x$, there exists a time $t > T$ such that $\Phi_t(x) \in V_x$. The vector field $X \in \mathcal{X}(M)$ is called* positively Poisson stable *if the set of positively Poisson stable points of $X$ is dense in $M$.*

**Definition 3.4.** *A point $x \in M$ is called a* nonwandering point *of $X \in \mathcal{X}(M)$ if for any $T > 0$ and for any neighborhood $V_x$ of $x$, there exists a time $t > T$ such that $\Phi_t(V_x) \cap V_x \neq \emptyset$.*

Let $\Gamma_X$ be the nonwandering set of $X$, i.e., the set of all nonwandering points of $X$. Then we have the following result [18].

**Theorem 3.5.** *The nonwandering set of a positively Poisson stable vector field $X$ is the entire $M$, that is, $\Gamma_X = M$.*

*Proof.* Let $x \in M$ be given. We want to prove that for any neighborhood $V_x$ of $x$ and for any $T > 0$, there exists a time $t > T$ such that $\Phi_t(V_x) \cap V_x \neq \emptyset$.

Let $S_X$ denote the set of positively Poisson stable points for $X \in \mathcal{X}(M)$. By definition, the closure of $S_X$ equals to $M$, i.e., $\overline{S_X} = M$. This implies that there exists a positively Poisson stable point $y$ in $V_x$. Now $V_x$ is also a neighborhood for $y$ and because $y$ is positively Poisson stable we have that for all $T > 0$, there exists a time $t > T$ such that $\Phi_t(y) \in V_x$. Hence $\Phi_t(V_x) \cap V_x \neq \emptyset$. Thus $x$ is a nonwandering point of $X$. Since $x$ was arbitrarily chosen, it follows that $\Gamma_X = M$, as required. □

Positive Poisson stability of a vector field is hence a sufficient condition for the nonwandering set to be the entire manifold. Since the converse is not true one introduces a weaker definition.

**Definition 3.6.** *A vector field is called weakly positively Poisson stable (WPPS) if its nonwandering set equals M, i.e., $\Gamma_X = M$.*

A natural question is when a vector field on a manifold is weakly positively Poisson stable (WPPS). For this purpose we recall two classical theorems. The first one is due to Liouville and states that the flow of a Hamiltonian vector field necessarily preserves the phase volume defined by the $n$th power of the symplectic form on $M$, where $2n = \dim M$. The second theorem is due to Poincaré. Let $\Omega$ be a volume form on $M$ and denote by $m_\Omega$ the associated Borel measure on an arbitrary manifold $M$. Let $X$ be a complete (time independent) vector field on $M$ whose flow $\{\Phi_t\}_{t \in \mathbb{R}}$ preserves the volume, i.e., $\Phi_t^* \Omega = \Omega$ for all $t \in \mathbb{R}$. Suppose $A$ is a measurable subset of $M$ with $0 < m_\Omega(A) < \infty$ which is also invariant under $\{\Phi_t\}_{t \in \mathbb{R}}$, i.e., $\Phi_t(A) \subset A$. The *Poincaré Recurrence Theorem* states that for each measurable subset $B$ of $A$ with $m_\Omega(B) > 0$ and for any $T > 0$, there exists some $t > T$ such that $\Phi_t(B) \cap B \neq \emptyset$. For the proof of these theorems see, for example [1], [2], [4], [19].

An immediate consequence is the following proposition.

**Proposition 3.7.** *Let $(M, \Omega)$ be a compact manifold with a volume form $\Omega$ and let $X$ be a time-independent vector field such that its flow preserves the volume form. Then it is a (WPPS) vector field.*

The link between the (WPPS) condition and controllability is given by the following theorem which is due to Kuang-You Lian, Li-Sheng Wang, and Li-Chen Fu [18]. Earlier versions of this theorem, where the hypothesis required $f$ to be Poisson stable, are due to Lorby [20] and Bonnard [8].

**Theorem 3.8.** *Suppose that $f$ is a (WPPS) vector field. The system (3.1) is controllable if and only if (LARC) holds.*

Now we give the setting for our result. Let $G$ be a Lie group acting freely properly and symplectically on a $2n$-dimensional connected symplectic manifold $(M, \omega)$. Suppose that the action admits an associated $Ad^*$-equivariant momentum map $J : M \longrightarrow \mathfrak{g}^*$. Consider on $M$ the nonlinear control system

$$\dot{x} = X_H(x) + \sum_{i=1}^{m} g_i(x) u_i \tag{3.2}$$

where $X_H$ is an Hamiltonian vector field, $g_1, ..., g_m \in \mathcal{X}(M)$ are $G$-invariant smooth vector fields and the control $u := (u_1, ..., u_m) : (0, \infty) \longrightarrow B \subset \mathbb{R}^m$ is a measurable function with values in a bounded subset $B$ of $\mathbb{R}^m$. Then the system (3.2) will naturally induce the nonlinear control system on $\left(M/G, \{\cdot, \cdot\}_{M/G}\right)$

$$\dot{\widetilde{x}} = X_{\widetilde{H}}(\widetilde{x}) + \sum_{i=1}^{m} \widetilde{g}_i(\widetilde{x}) u_i, \tag{3.3}$$

where $X_{\widetilde{H}}$ is a Hamiltonian vector field with respect to the Poisson bracket $\{\cdot, \cdot\}_{M/G}$ and Hamiltonian $\widetilde{H}$ given by $H = \widetilde{H} \circ \pi$ where, $\pi : M \longrightarrow M/G$ is the canonical projection.

**Theorem 3.9.** *Suppose that the system (3.3) verifies (LARC).*

*(a) If the momentum map $J : M \longrightarrow \mathfrak{g}^*$ is a proper map, then the system (3.3) is controllable.*

*(b) If the momentum map is not proper but the Lie group $G$ is compact and if there exists a continuous proper map $V : M/G \longrightarrow \mathbb{R}$ which is constant along the trajectories of $X_{\widetilde{H}}$, then the system (3.3) is controllable.*

*Proof.* The strategy to prove the controllability of (3.3) is to show that $X_{\widetilde{H}}$ is (WPPS) and then the conclusion follows from Theorem 3.8.

(a) As subsets of $M/G$, the symplectic leaves are $M_{\mathcal{O}_\mu}$ or, equivalently, $j_\mu(M_\mu)$ and the symplectic form is given by $\omega_\mu^{\maltese}$. Because $J$ is a proper map, the set $J^{-1}(\mu)$ is a compact submanifold of $M$. Thus $M_\mu$ is a compact manifold which implies that the injective immersion $j_\mu$ is in fact an embedding. So, the immersed topology and the induced topology on $M_{\mathcal{O}_\mu}$ coincide and therefore the symplectic leaves are compact submanifolds of $M/G$.

The flow of $X_{\widetilde{H}}$ leaves the symplectic leaves of $M/G$ invariant. On each such leaf the flow of $X_{\widetilde{H}}$ preserves the symplectic form $\omega_\mu^{\maltese}$ and so, by Liouville's theorem, it also preserves the induced volume form.

Let $\widetilde{x}_0$ be an arbitrary point in $M/G$. It belongs to a symplectic leaf denoted by $L_{\mu_0}$. By Proposition 3.7, $\widetilde{x}_0$ is a nonwandering point for the Hamiltonian vector field given by the restriction of $X_{\widetilde{H}}$ to the symplectic leaf $L_{\mu_0}$. This implies that $\widetilde{x}_0$ is also a nonwandering point for $X_{\widetilde{H}}$. As $\widetilde{x}_0$ was chosen arbitrarily, we obtain that $X_{\widetilde{H}}$ is (WPPS).

(b) For compact $G$, the coadjoint orbits are submanifolds of $\mathfrak{g}^*$ and $J$ is transversal to the coadjoint orbits that lie in its image (since by hypothesis, the action is free). So $J^{-1}(\mathcal{O}_\mu)$ is a submanifold of $M$ and by Theorem 2.2(i) the immersed topology and the induced topology on $M_{\mathcal{O}_\mu}$ coincide.

As before, let $\widetilde{x}_0$ be an arbitrary point in $M/G$. It belongs to a symplectic leaf denoted by $L_{\mu_0}$. Combining the above observation with the result of Proposition 2.1 we conclude that $L_{\mu_0}$ is a closed submanifold of $M/G$. If $\dim L_{\mu_0} = 0$ then, because $L_{\mu_0}$ is connected, it follows that it equals the single point set $\{\widetilde{x}_0\}$. Thus $\widetilde{x}_0$ is an equilibrium of $X_{\widetilde{H}}$ and so it is a nonwandering point.

If $\dim L_{\mu_0} > 0$, let $m_{\mu_0}$ be the Borel measure associated to the Liouville volume form induced by the symplectic form $\omega_{\mu_0}^{\maltese}$. Let $c_0 := V(\widetilde{x}_0)$. If $V(L_{\mu_0}) = c_0$ then $L_{\mu_0} \subset V^{-1}(c_0)$ and because $V^{-1}(c_0)$ is a compact subset and $L_{\mu_0}$ is closed in $M/G$, we conclude that $L_{\mu_0}$ is compact. By Proposition 3.7, $\widetilde{x}_0$ is a nonwandering point.

The other possibility is that the image of $V_{|L_{\mu_0}}$ is a non-degenerate connected interval $I \subset \mathbb{R}$ that contains $c_0$. Taking a small interval around $c_0$, $[-\epsilon + c_0, c_0 + \epsilon]$ we have that the flow invariant set $K = L_{\mu_0} \cap V^{-1}([-\epsilon + c_0, c_0 + \epsilon])$ is a compact subset of $M/G$ because $L_{\mu_0}$ is closed and $V^{-1}([-\epsilon + c_0, c_0 + \epsilon])$ is compact in $M/G$. This shows that $m_{\mu_o}(K) < \infty$. Also, $K$ contains an open subset of $L_{\mu_0}$, for example $V_{|L_{\mu_0}}^{-1}((-\epsilon + c_0, c_0 + \epsilon))$, implying that $0 < m_{\mu_o}(K)$.

Consider now $U_{\widetilde{x}_0}$ an arbitrary open neighborhood of $\widetilde{x}_0$ in $M/G$. The set $K \cap U_{\widetilde{x}_0}$ is a Borel subset of $L_{\mu_0}$ and because it contains a non-empty open subset of $L_{\mu_0}$, for example $U_{\widetilde{x}_0} \cap V_{|L_{\mu_0}}^{-1}((-\epsilon + c_0, c_0 + \epsilon))$, it follows that $m_{\mu_o}(K \cap U_{\widetilde{x}_0}) > 0$. By the Poincaré Recurrence Theorem, for all $T > 0$ there exists $t > T$ such that $\Phi_t \left( K \cap U_{\widetilde{x}_0} \right) \cap (K \cap U_{\widetilde{x}_0}) \neq \emptyset$ which, in particular, implies that $\Phi_t \left( U_{\widetilde{x}_0} \right) \cap U_{\widetilde{x}_0} \neq \emptyset$. Consequently, $\widetilde{x}_0$ is a nonwandering point for $X_{\widetilde{H}}$.

As before, since $\widetilde{x}_0$ was arbitrary in $M/G$ it follows that $X_{\widetilde{H}}$ is (WPPS). □

**Remark 3.10.** Observe that for the controllability of (3.3) it is not necessary for the vector fields $\widetilde{g}_i \in \mathcal{X}(M/G)$ to be induced by some $G$-invariant vector fields on $M$.

**Remark 3.11.** Under the conditions of the previous theorem, part (b), the system given by (3.3) is controllable. If, in addition, the system (3.2) is accessible then using a theorem proved by San Martin and Crouch [29] it follows that (3.2) is also controllable. If the hypotheses of the theorem in San Martin and Crouch [29] do not hold, then a possible starting point for studying the relationship between controllability of the original and reduced systems is given in Grizzle and Marcus [11].

## 4. Examples

We will illustrate the theory with three examples. In all of them we will use the following well known lemmas to prove the properness of the integrals of motion.

**Lemma 4.1.** *Let $f : \mathbb{R}^n \to \mathbb{R}^k$ be a continuous function. Then $f$ is proper if and only if*

$$\lim_{\|x\| \to \infty} \|f(x)\| = +\infty.$$

*Proof.* Suppose that $f$ is proper. If $\lim_{\|x\| \to \infty} \|f(x)\| \neq +\infty$, there exists a sequence $\{x_n\}_{n \in \mathbb{N}}$ and a constant $M > 0$ such that $\|x_n\| \to \infty$ and $\|f(x_n)\| \leq M$. Thus $\{x_n\}_{n \in \mathbb{N}}$ lies in the inverse image by $f$ of the closed ball of radius $M$ which is a compact set in $\mathbb{R}^n$ because $f$ is assumed to be proper. Hence $\{x_n\}_{n \in \mathbb{N}}$ contains a convergent subsequence. However, $\|x_n\| \to \infty$ which is a contradiction.

Conversely, assume that $\lim_{\|x\| \to \infty} \|f(x)\| = +\infty$ and let $K \subset \mathbb{R}^k$ be a compact subset. The set $f^{-1}(K)$ is closed since $f$ is continuous. To conclude that $f^{-1}(K)$ is compact we shall show that it is also bounded. If not, there would exist a sequence $\{x_n\}_{n \in \mathbb{N}} \subset f^{-1}(K)$ such that $\|x_n\| \to \infty$. By hypothesis, $\|f(x_n)\| \to \infty$, which contradicts the fact that $f(x_n) \in K$ which is bounded. □

**Lemma 4.2.** *Let $M, N, P$ be Hausdorff topological spaces and let $f : M \to N$ and $g : N \to P$ be two continuous functions. If $g \circ f : M \to P$ is proper, then $f$ is also proper.*

*Proof.* Let $K \subset M$ be a compact subset. Then $g(K)$ is compact in $P$ and hence $(g \circ f)^{-1}(g(K))$ is compact in $M$. Since $f^{-1}(K) \subset (g \circ f)^{-1}(g(K))$ is closed, it follows that it is also compact. □

**Example 1**. The motion of three point vortices for an ideal inviscid incompressible fluid in the plane is given by the equations

(4.1)
$$\begin{cases} \dot{x}_j = -\frac{1}{2\pi} \sum_{\substack{i=1 \\ i \neq j}}^{3} \Gamma_i (y_j - y_i)/r_{ij}^2 \\ \dot{y}_j = -\frac{1}{2\pi} \sum_{\substack{i=1 \\ i \neq j}}^{3} \Gamma_i (x_j - x_i)/r_{ij}^2, \end{cases}$$

where $r_{ij}^2 = (x_i - x_j)^2 + (y_i - y_j)^2$ and $\Gamma_1, \Gamma_2, \Gamma_3$ are non-zero constants, the circulations given by the corresponding point vortices. Kirchoff [16] noted that (4.1) can be written in the form

$$\Gamma_j \frac{dx_j}{dt} = \frac{\partial H}{\partial y_j}$$
$$\Gamma_j \frac{dy_j}{dt} = -\frac{\partial H}{\partial x_j},$$

where

$$H(x_1, x_2, x_3, y_1, y_2, y_3) = -\frac{1}{4\pi} \sum_{\substack{i=1 \\ i \neq j}}^{3} \Gamma_i \Gamma_j \log r_{ij}$$

is the Hamiltonian and the symplectic form is given by

(4.2)
$$\Omega = \sum_{i=1}^{3} \Gamma_i dx_i \wedge dy_i.$$

Consider the diagonal action of $SE(2)$ on $(\mathbb{R}^2)^3$ whose associated momentum map $J : \mathbb{R}^6 \to \mathbb{R}^3$ relative to the symplectic form (4.2) is given by

$$J(\mathbf{x}, \mathbf{y}) = \left( -\frac{1}{2} \sum_{i=1}^{3} \Gamma_i (x_i^2 + y_i^2), \sum_{i=1}^{3} \Gamma_i y_i, -\sum_{i=1}^{3} \Gamma_i x_i \right).$$

We will identify $\mathbb{R}^2$ with $\mathbb{C}$ and define

$$Q := \{(u, 0) \mid u \in \mathbb{C}\} \cup \{(0, v) \mid v \in \mathbb{C}\} \cup \{(w, w) \mid w \in \mathbb{C}\}.$$

On the set $S = \mathbb{C}^3 \setminus Q$ the action of $SE(2)$ is free and proper. If the vortex strengths $\Gamma_1, \Gamma_2, \Gamma_3$ have the same signs then by applying Lemma 4.2 and Lemma 4.1 it follows that $J$ is a proper map and we are in the first case of Theorem 3.9.

The matrix of the Poisson bracket on the reduced space $S/SE(2) \approx T := \mathbb{R}^3 \setminus (\{(0, 0, c) \mid c \in \mathbb{R}\} \cup \{(a, 0, 0) \mid a \geq 0\})$ is

$$4 \begin{bmatrix} 0 & 2a_3 & -2a_2 \\ -2a_3 & 0 & 2a_1 - \|\mathbf{a}\| \\ 2a_2 & -2a_1 + \|\mathbf{a}\| & 0 \end{bmatrix}$$

and the reduced Hamiltonian is

$$\widetilde{H}(a_1, a_2, a_3) = -\frac{1}{4\pi} (\Gamma_1 \Gamma_2 \log((a_3 + \|\mathbf{a}\|)/2) + \Gamma_1 \Gamma_3 \log((-a_3 + \|\mathbf{a}\|)/2)$$
$$+ \Gamma_2 \Gamma_3 \log(-a_1 + \|\mathbf{a}\|)).$$

The reduced equations are

$$\dot{a}_1 = \frac{2}{\pi}\left(\Gamma_1\Gamma_2\frac{a_2}{(a_3+\|\mathbf{a}\|)} - \Gamma_1\Gamma_3\frac{a_2}{(-a_3+\|\mathbf{a}\|)}\right)$$

$$\dot{a}_2 = \frac{1}{\pi}\left(\Gamma_1\Gamma_2\frac{(-2a_1+a_3+\|\mathbf{a}\|)}{(a_3+\|\mathbf{a}\|)} + \Gamma_1\Gamma_3\frac{(2a_1+a_3-\|\mathbf{a}\|)}{(-a_3+\|\mathbf{a}\|)} + \Gamma_2\Gamma_3\frac{a_3}{(-a_1+\|\mathbf{a}\|)}\right)$$

$$\dot{a}_3 = \frac{1}{\pi}\left(\Gamma_2\Gamma_3\frac{a_2}{(-a_1+\|\mathbf{a}\|)} - \Gamma_1\Gamma_2\frac{a_2}{(a_3+\|\mathbf{a}\|)} - \Gamma_1\Gamma_3\frac{a_2}{(-a_3+\|\mathbf{a}\|)}\right)$$

and consider the reduced controlled system

$$\dot{a}_1 = \frac{2}{\pi}\left(\Gamma_1\Gamma_2\frac{a_2}{(a_3+\|\mathbf{a}\|)} - \Gamma_1\Gamma_3\frac{a_2}{(-a_3+\|\mathbf{a}\|)}\right) + u_1$$

$$\dot{a}_2 = \frac{1}{\pi}\left(\Gamma_1\Gamma_2\frac{(-2a_1+a_3+\|\mathbf{a}\|)}{(a_3+\|\mathbf{a}\|)} + \Gamma_1\Gamma_3\frac{(2a_1+a_3-\|\mathbf{a}\|)}{(-a_3+\|\mathbf{a}\|)} + \Gamma_2\Gamma_3\frac{a_3}{(-a_1+\|\mathbf{a}\|)}\right) + u_2$$

$$\dot{a}_3 = \frac{1}{\pi}\left(\Gamma_2\Gamma_3\frac{a_2}{(-a_1+\|\mathbf{a}\|)} - \Gamma_1\Gamma_2\frac{a_2}{(a_3+\|\mathbf{a}\|)} - \Gamma_1\Gamma_3\frac{a_2}{(-a_3+\|\mathbf{a}\|)}\right) + u_3,$$

where the control $u := (u_1, u_2, u_3) : (0, \infty) \to B \subset \mathbb{R}^3$ is a measurable function with values in a bounded subset $B \subset \mathbb{R}^3$.

It is easy to check that the vector fields $X_{\widetilde{H}}, \frac{\partial}{\partial a_1}, \frac{\partial}{\partial a_2}, \frac{\partial}{\partial a_3}$ verify (LARC) and as a result of Theorem 3.9 (a), we conclude that the reduced system above is controllable.

**Example 2**. The next example will be the resonant three-wave interaction. This is a Hamiltonian system with the phase space $P \equiv \mathbb{R}^6 = \mathbb{C}^3$, equipped with the symplectic structure

$$\omega = \sum_{j=1}^{3} \frac{1}{s_j \gamma_j} dq_j \wedge dp_j,$$

where $s_1, s_2, s_3 \in \{-1, 1\}$ and $\gamma_1, \gamma_2, \gamma_3 \in \mathbb{R}$ are parameters subject to the constraint $\gamma_1 + \gamma_2 + \gamma_3 = 0$. We will restrict our attention to the particular case where $(s_1, s_2, s_3) = (1, 1, 1)$ and $(\gamma_1, \gamma_2, \gamma_3) = (1, 1, -2)$.

In standard coordinates on $\mathbb{C}^3$ the Hamiltonian is given by

$$H(z_1, z_2, z_3) = -\frac{1}{2}(\bar{z}_1 z_2 \bar{z}_3 + z_1 \bar{z}_2 z_3).$$

This Hamiltonian is invariant under the action of the compact Lie group $G \equiv S^1 \times S^1$ on $P$ given by

$$(e^{i\theta_1}, e^{i\theta_2}) \cdot (z_1, z_2, z_3) = (e^{-i\theta_1} z_1, e^{-i(\theta_1+\theta_2)} z_2, e^{-i\theta_2} z_2), \qquad 0 \leqslant \theta_j < 2\pi.$$

The momentum map for this action is $J : P \to \mathfrak{g}^* \cong \mathbb{R}^2$

$$J(z_1, z_2, z_3) = \left(\frac{1}{2}(|z_1|^2 + |z_2|^2), \frac{1}{2}(|z_1|^2 - \frac{1}{2}|z_2|^2)\right).$$

In order to have a free action we have to remove certain points from the phase space $P$ and replace it with $P \equiv \mathbb{C}\backslash\{0\} \times \mathbb{C} \times \mathbb{C}\backslash\{0\}$.

The matrix of the Poisson bracket on the reduced space $P/G \approx Q = \mathbb{R}^2 \times (0, \infty)^2$ is

$$\begin{bmatrix} 0 & 1 & -\frac{p}{a} & 2\frac{p}{b} \\ -1 & 0 & \frac{q}{a} & -2\frac{q}{b} \\ \frac{p}{a} & -\frac{q}{a} & 0 & 0 \\ -2\frac{p}{b} & 2\frac{q}{b} & 0 & 0 \end{bmatrix}$$

and the reduced Hamiltonian is

$$\widetilde{H}(q, p, a, b) = -abq.$$

The reduced equations of motion are

(4.3)
$$\begin{aligned} \dot{q} &= \frac{qpb}{a} - 2\frac{qpa}{b} \\ \dot{p} &= -ab - \frac{q^2 b}{a} + 2\frac{q^2 a}{b} \\ \dot{a} &= -pb \\ \dot{b} &= 2pa. \end{aligned}$$

A constant of motion for the system (4.3) is given by the function $V : Q \to \mathbb{R}$, $V(q, p, a, b) = q^2 + p^2 + a^2 + b^2$ which is proper by Lemma 4.1.

Consider now the following reduced controlled system

(4.4)
$$\begin{aligned} \dot{q} &= \frac{qpb}{a} - 2\frac{qpa}{b} + u_1 \\ \dot{p} &= -ab - \frac{q^2 b}{a} + 2\frac{q^2 a}{b} + u_2 \\ \dot{a} &= -pb \\ \dot{b} &= 2pa + u_3, \end{aligned}$$

where the control $u := (u_1, u_2, u_3) : (0, \infty) \to B \subset \mathbb{R}^3$ is a measurable function with values in a bounded subset $B$. A short computation shows that the vector fields $\{\partial/\partial q, \partial/\partial p, \partial/\partial b, [\partial/\partial b, [\partial/\partial p, X_{\widetilde{H}}]]\}$ at every point $x \in Q$ generate the tangent space $T_x Q$, which proves that the system (4.4) verifies (LARC). Applying now the result of Theorem 3.9 (b), we obtain that (4.4) is controllable.

**Example 3**. We will study the controllability of the reduced system of two coupled planar rigid bodies. We take the description of the system given in [32]. After the reduction to the center of mass frame we have the configuration space $S^1 \times S^1$ with the diagonal action of $S^1$. The phase space is $T^*(S^1 \times S^1)$ with the canonical symplectic form of a cotangent bundle. The momentum map for the lifted action of $S^1$ is given by

$$J((\theta_1, \mu_1), (\theta_2, \mu_2)) = \mu_1 + \mu_2.$$

Following Krishnaprasad and Marsden [17] the reduced Poisson space is

$$P := T^*(S^1 \times S^1)/S^1 \cong S^1 \times \mathbb{R}^2$$

and if we chose coordinates $(\theta, \mu_1, \mu_2)$ on $P$ the matrix of the Poisson bracket is given by

$$\begin{bmatrix} 0 & -1 & 1 \\ 1 & 0 & 0 \\ -1 & 0 & 0 \end{bmatrix}.$$

The reduced Hamiltonian is given by the formula

$$H = \frac{1}{2\triangle}(\widetilde{I_2}\mu_1^2 - 2\varepsilon\lambda(\theta)\mu_1\mu_2 + \widetilde{I_1}\mu_2^2),$$

where $\triangle = \widetilde{I_1}\widetilde{I_2} - \varepsilon^2(\lambda(\theta))^2 > 0$ and

| | |
|---|---|
| $d_i$ | distance from the hinge to the center of mass of body $i = 1, 2$ |
| $\theta$ | joint angle from body 1 to body 2 |
| $\lambda(\theta)$ | $d_1 d_2 \cos\theta$ |
| $m_i$ | mass of body $i = 1, 2$ |
| $\varepsilon$ | $m_1 m_2/(m_1 + m_2)$ = reduced mass |
| $I_i$ | moment of inertia of body $i$ about its center of mass |
| $\widetilde{I_i}$ | $I_i + \varepsilon d_i^2$, $i = 1, 2$ (augmented moments of inertia). |

To apply Theorem 3.9 we need to show that $H$ is a proper function. To do this, we need the following lemma.

**Lemma 4.3.** *Let $f : K \to \mathbb{R}$ and $g : \mathbb{R}^n \to \mathbb{R}$ be two continuous functions, where $K$ is compact and $g$ is a proper function. Then the function $h : K \times \mathbb{R}^n \to \mathbb{R}$ given by $h(x, y) := f(x)g(y)$ is a proper function.*

*Proof.* We shall prove that $h^{-1}([a, b])$ is compact in $K \times \mathbb{R}^n$. Let $z_n := (x_n, y_n)$ be an arbitrary sequence in $h^{-1}([a, b])$. Since $K$ is compact, we can assume that $\{x_n\}_{n \in \mathbb{N}}$ is convergent. Because $\{f(x_n)g(y_n)\}_{n \in \mathbb{N}} \subset [a, b]$ and $\{f(x_n)\}_{n \in \mathbb{N}}$ is bounded, the sequence $\{g(y_n)\}_{n \in \mathbb{N}}$ is also bounded and hence there are $a', b' \in \mathbb{R}$ such that $\{g(y_n)\}_{n \in \mathbb{N}} \subset [a', b']$. Therefore $\{y_n\}_{n \in \mathbb{N}} \subset g^{-1}([a', b'])$ which is a compact set in $\mathbb{R}^n$ because $g$ is a proper function. Consequently there is a convergent subsequence of $\{y_n\}_{n \in \mathbb{N}}$. The corresponding subsequence of $\{z_n\}_{n \in \mathbb{N}}$ is convergent which proves that $h^{-1}([a, b])$ is compact. $\square$

To apply this lemma we write $H$ in the form

$$H = \frac{1}{2\triangle}\left(\left(\sqrt{\widetilde{I_2}}\,\mu_1 - \frac{\varepsilon\lambda(\theta)}{\sqrt{\widetilde{I_2}}}\mu_2\right)^2 + \left(\widetilde{I_1} - \frac{\varepsilon^2\lambda^2(\theta)}{\widetilde{I_2}}\right)\mu_2^2\right).$$

Since

$$\widetilde{I_1} - \frac{\varepsilon^2\lambda^2(\theta)}{\widetilde{I_2}} > 0$$

the smooth change of variables $(\theta, \mu_1, \mu_2) \mapsto (\theta, X, Y)$, where

$$X := \sqrt{\widetilde{I_2}}\, \mu_1 - \frac{\varepsilon \lambda(\theta)}{\sqrt{\widetilde{I_2}}} \mu_2$$

$$Y := \left( \widetilde{I_1} - \frac{\varepsilon^2 \lambda^2(\theta)}{\widetilde{I_2}} \right)^{1/2} \mu_2$$

transforms $H$ to the function $\frac{1}{2\triangle} \left( X^2 + Y^2 \right)$ with $1/2\triangle$ defined on $S^1$. This function is proper by Lemma 4.1 and Lemma 4.3. Thus $H$ is a proper integral of motion for the reduced system.

Now we consider the following reduced controlled system

$$\begin{aligned}
\dot{\theta} &= -\frac{\partial H}{\partial \mu_1} + \frac{\partial H}{\partial \mu_2} + u_1 \\
\dot{\mu_1} &= \frac{\partial H}{\partial \theta} + u_2 \\
\dot{\mu_2} &= -\frac{\partial H}{\partial \theta} + u_3
\end{aligned}$$

where the control $u := (u_1, u_2, u_3) : (0, \infty) \to B \subset \mathbb{R}^3$ is a measurable function with values in a bounded subset $B$. It is easy to see that the vector fields $X_H, \frac{\partial}{\partial \theta}, \frac{\partial}{\partial \mu_1}, \frac{\partial}{\partial \mu_2}$ verify (LARC) and as a result of Theorem 3.9 (b), we obtain that the reduced controlled system above is controllable.

**Acknowledgment**. We would like to thank Vlad Timofte and Ioan Casu for carefully reading the manuscript and valuable comments which helped us to improve the exposition. This research was partially supported by the European Commission and the Swiss Federal Government through funding for the Research Training Network *Mechanics and Symmetry in Europe* (MASIE) (T.S.R.) as well as the Swiss National Science Foundation through FNS Grant 20-61228.00 (P.B. and T.S.R.).